# Безградиентные прокс-методы с неточным оракулом для негладких задач выпуклой стохастической оптимизации на симплексе[1]

*Гасников А.В. (к.ф.-м.н., ИППИ РАН, ПреМоЛаб ФУПМ МФТИ) gasnikov@yandex.ru*

*Лагуновская А.А. (ИПМ им.М.В.Келдыша РАН, МФТИ) a.lagunovskaya@phystech.edu*

*Усманова И.Н. (ПреМоЛаб ФУПМ МФТИ) ilnura94@gmail.com*

*Федоренко Ф.А. (Кафедра МОУ ФУПМ МФТИ) f.a.fedorenko@gmail.com*

**Аннотация**

В работе предложена безградиентная модификация метода зеркального спуска решения задач негладкой стохастической выпуклой оптимизации на единичном симплексе. Особенностью постановки является допущение, что реализации значений функции нам доступны с небольшими шумами. Цель данной работы – установить скорость сходимости предложенного метода, и определить, при каком уровне шума, факт его наличия не будет существенно сказываться на скорости сходимости.

**Ключевые слова:** метод зеркального спуска, безградиентные методы, методы с неточным оракулом, стохастическая оптимизация.

---

## 1. Введение

Представим себе, следуя Ю.Е. Нестерову [1], что некоторый человек может характеризовать свое состояние вектором

$$x \in S_n(1) = \left\{ x \geq 0: \quad \sum_{i=1}^{n} x_i = 1 \right\}.$$

Насколько это состояние хорошее он может оценить, посчитав значение своей функции потерь $f(x)$ на этом векторе. К сожалению, рассматриваемый человек существенно ограничен в своих возможностях, поэтому посчитать (суб-)градиент этой функции он не может. Более того, значение функции он может посчитать лишь с неконтролируемым им шумом уровня $\delta$. Функция потерь (заданная в некоторой окрестности симплекса) предполагается выпуклой, но необязательно гладкой, с равномерно ограниченной нормой субградиента $\|\nabla f(x)\|_\infty \leq M$. Человек стремится оказаться в состоянии с наименьшими потерями $f_* = \min_{x \in S_n(1)} f(x)$, действуя итерационно по следующему простому правилу:

- выбрать случайно направление;
- сдвинуться с некоторым (небольшим) шагом из текущего состояния по этому направлению;
- посчитать значение функции в новом состоянии;
- исходя из значения функции в этих двух состояниях (состояния, в котором находились, и состояния, в котором оказались), определить новое состояние (из симплекса), в которое следует перейти.

В связи с описанным естественным способом действий возникает ряд вопросов. Например, как выбирать новое состояние, с целью минимизации числа обращений к «оракулу» (см. раздел 3) за значением функции? Как именно "случайно" стоит выбирать направление? Как скажется зашумленность выдаваемых значений функции на это число обращений? Можно ли приблизиться к нижним оценкам требуемого числа обращений для достижения $f_*$ с точностью (по функции) $\varepsilon$?

В статье мы рассматриваем еще более общую постановку, когда оракул может выдавать не значение функции, а лишь несмещенную (или не сильно смещенную, смещение контролируется уровнем шума $\delta$) оценку этого значения $f(x;\eta)$:

$$E_\eta [f(x;\eta)] = f(x).$$

В такой общности мы постараемся ответить на сформулированные вопросы. В частности, будет предложена процедура, требующая в случае гладкой функции $f(x)$ (точнее липшицевости градиента $f(x;\eta)$ по $x$)

$$\mathrm{O}\left(\frac{M^2 n \ln n}{\varepsilon^2}\right)$$

обращений к оракулу за реализацией функции $f(x;\eta)$, что с точностью до логарифмического множителя соответствует нижней оценке [2]. Предложенный в работе подход, позволяет также при некоторых дополнительных предположениях заметно улучшить приведенную оценку.

Отметим, что если бы вместо значения функции оракул выдавал стохастический градиент или хотя бы стохастическую производную функции по направлению, то ответы были бы, соответственно:

$$\mathrm{O}\left(\frac{M^2 \ln n}{\varepsilon^2}\right), \mathrm{O}\left(\frac{M^2 n \ln n}{\varepsilon^2}\right),$$

что также с точностью до логарифмических множителей соответствует нижним оценкам (см., например, [2, 3]).

В целом проблематика работы восходит к статье [4] (см. также [5, 6, 7, 8]). В п. 2 мы описываем известные результаты о сходимости метода зеркального спуска (МЗС) для задач стохастической оптимизации [7], которые нам понадобятся в дальнейшем. В п. 3 мы вводим неточный оракул, выдающий зашумленное значения реализации функции. Исходя из такой (частичной) информации в п. 3 предлагаются различные рандомизированные (безградиентные) обобщения МЗС. Рандомизация заключается в выборе случайного направления и вычислении (с помощью оракула) вместо стохастического субградиента стохастической дискретной производной функции по этому направлению [4]. Основные степени свободы, на которых можно играть: способ выбора случайного направления (в работе обсуждаются равномерное распределение на евклидовом шаре и равномерные распределения на шарах в $l_1$ и $l_\infty$ нормах) и выбор шага дискретизации. В отсутствии шума выгоднее всего этот шаг стремить к нулю, т.е. просто вычислять стохастическую производную по направлению. Однако мы допускаем шум, и хотим понять, при каком максимально допустимом уровне шума оценки сохранят свой вид, скажем, в таких категориях: число итераций возрастет не более чем в два раза. В пп. 2, 3 гладкость не предполагается. В п. 4 на примере изучения стохастических спусков по случайным направлением демонстрируется увеличение скорости сходимости, связанное с наличием гладкости. В п. 5 результаты п. 4 переносятся на стохастические безградиентные методы, т.е. по сути, на гладкий вариант постановки задачи из п. 3. Наличие гладкости дает ускорение в пп. 4, 5 приблизительно в $n$ раз.

## 2. Метод зеркального спуска для задач стохастической оптимизации с неточным оракулом

Рассмотрим задачу стохастической оптимизации

(1) $$f(x) = E_\eta\left[f(x;\eta)\right] \to \min_{x \in S_n(1)}.$$

Здесь $\eta$ – случайная величина, $E_\eta\left[f(x;\eta)\right]$ – математическое ожидание "взятое по $\eta$", то есть при фиксированном $x$, при этом далее допускается, что в такой записи $x$ может быть случайным вектором. В таком случае математическое ожидание берётся только по $\eta$ (случайность в $x$ "фиксируется"). Если математическое ожидание берётся по $x$ (первое неравенство в теореме 1), то нижний индекс $\eta$ опускаем.

Обозначим

$$f_* = \min_{x \in S_n(1)} f(x) = \min_{x \in S_n(1)} E_\eta\left[f(x;\eta)\right].$$

**Замечание 1.** Везде далее мы будем использовать обозначения обычного градиента для субградиента. Запись $\nabla_x f(x;\eta)$ в вычислительном контексте (например, в итерационной процедуре (2) ниже) означает какой-либо измеримый селектор стохастического субдифференциала [9], а если в контексте проверки условий (например, в условии 2 или условии 3 ниже), то $\nabla_x f(x;\eta)$ пробегает все элементы стохастического субдифференциала.

Для формулировки основной теоремы этого пункта нам понадобятся следующие **условия**:

1. $f(x;\eta)$ – выпуклая функция по $x$ (в действительности, с некоторыми оговорками [9], достаточно только выпуклости $f(x)$);

2. Стохастический субградиент $\nabla_x f(x;\eta)$ [9] удовлетворяет условию (тождественно по $x$):

$$E_\eta\left[\nabla_x f(x;\eta)\right] \equiv \nabla_x E_\eta\left[f(x;\eta)\right];$$

3. $\left\|\nabla_x f(x;\eta)\right\|_\infty \leq M$ – равномерно, с вероятностью 1.

Для справедливости части утверждений достаточно требовать вместо условия 3 одно из следующих (более слабых) условий:

$$\text{а) } E_\eta\left[\|\nabla_x f(x;\eta)\|_\infty^2\right] \leq M^2; \quad \text{б) } E_\eta\left[\exp\left(\frac{\|\nabla_x f(x;\eta)\|_\infty^2}{M^2}\right)\right] \leq \exp(1).$$

Для решения задачи (1) воспользуемся методом зеркального спуска (точнее двойственных усреднений) в форме [10, 11]. Положим $x_i^1 = 1/n$, $i = 1,...,n$. Пусть $t = 1,...,N-1$.

$$(2) \quad x_i^{t+1} = \frac{\exp\left(-\frac{1}{\beta_{t+1}}\sum_{k=1}^{t}\frac{\partial f(x^k;\eta^k)}{\partial x_i}\right)}{\sum_{l=1}^{n}\exp\left(-\frac{1}{\beta_{t+1}}\sum_{k=1}^{t}\frac{\partial f(x^k;\eta^k)}{\partial x_l}\right)}, \quad i = 1,...,n, \quad \beta_t = \frac{M\sqrt{t}}{\sqrt{\ln n}}.$$

Здесь $\{\eta^k\}$ – независимые одинаково распределенные (также как $\eta$) случайные величины.

Приводимая ниже теорема фактически установлена в работах [10, 11, 12]. Однако здесь мы непосредственно воспользовались формулировкой из работы [13].

**Теорема 1.** *Пусть справедливы условия 1, 2, 3.а, тогда*

$$E\left[f\left(\frac{1}{N}\sum_{k=1}^{N}x^k\right)\right] - f_* \leq \frac{1}{N}\sum_{k=1}^{N}E\left[f(x^k)\right] - f_* \leq 2M\sqrt{\frac{\ln n}{N}}.$$

*Пусть справедливы условия 1, 2, 3, тогда при $\Omega \geq 0$*

$$P_{x^1,...,x^N}\left\{\frac{1}{N}\sum_{k=1}^{N}f(x^k) - f_* \geq \frac{2M}{\sqrt{N}}\left(\sqrt{\ln n} + \sqrt{8\Omega}\right)\right\} \leq$$

$$\leq P_{x^1,...,x^N}\left\{f\left(\frac{1}{N}\sum_{k=1}^{N}x^k\right) - f_* \geq \frac{2M}{\sqrt{N}}\left(\sqrt{\ln n} + \sqrt{8\Omega}\right)\right\} \leq \exp(-\Omega).$$

**Замечание 2.** Если вместо условия 3 имеет место более слабое условие 3.б, то последняя формула останется верной, при небольшой корректировке:

$$\frac{2M}{\sqrt{N}}\left(\sqrt{\ln n} + \sqrt{8\Omega}\right) \to C\frac{M}{\sqrt{N}}\left(\sqrt{\ln n} + \Omega\right),$$

где константа $C \sim 10$. Приведенный результат можно обобщить и на более тяжелые хвосты [14].

## 3. Безградиентная модификация метода зеркального спуска для задач стохастической оптимизации с неточным оракулом

Введем понятие оракула, выдающего зашумленное значение функции $f(x)$, определенной[2] в $\mu_0$-окрестности $S_n(1)$.[3] При этом везде в дальнейшем под $f(x)$ и $f(x,\eta)$ мы будем понимать, соответственно, не зашумленные значение функции и ее несмещенной (тождественно по $x$) реализации. Наличие шума мы будем явно указывать, вводя его аддитивным образом.

**Предположение 1.** *Оракул выдает (на запрос, в котором указывается только одна точка $x$) $f(x,\eta)+\tilde{\delta}(\eta)$, где случайная величина $\eta$ независимо разыгрывается из одного и того же распределения, фигурирующего в постановке (1); случайная величина $\tilde{\delta}(\eta)$ (случайность может быть обусловлена не только зависимостью от $\eta$) не зависит от $x$ и ограничена по модулю известным нам числом $\delta$ – допустимым уровнем шума.*

Приведем одну из возможных мотивировок такого оракула. Предположим, что оракул может считать абсолютно точно значение (или реализацию) функции, но вынужден нам выдавать лишь конечное (предписанное) число первых бит (конечная мантисса). Таким образом, в последнем полученном бите есть некоторая неточность (причем неизвестно, по какому правилу оракул формирует этот последний выдаваемый значащий бит). Однако всегда можно прибавить (по mod 1) к этому биту случайно приготовленный (независимый) бит. В результате, не ограничивая общности, можно считать, что оракул последний бит выбирает просто случайно в независимости от отброшенного остатка.

**Предположение 2.** *В случае задач стохастической оптимизации принципиально важно, что разрешается на каждом шаге (итерации) обратиться к оракулу за значе-*

---

[2] Везде далее в статье мы будем предполагать, что $f(x)$ не просто определена в достаточно большой $\mu_0$-окрестности исходного множества, но и сохраняет все свои свойства в этой окрестности, в частности, выпуклость и константы Липшица.

[3] Все, что будет написано далее, можно перенести (без изменений итоговых формул с точностью до константного множителя) на случай более общего оракула, описанного в разделе 4 [14]. К сожалению, в [14] все равно относительно оракула делаются обременительные предположения. Впрочем, в этой же работе схематично показано, как можно распространить (с ужесточением условий на допустимый уровень шума) все, что далее будет написано на случай самого общего оракула, выдающего зашумленное значение функции (реализации функции). Об этом также написано в работе [15].

*ниями функции на одной реализации ($\eta$ одно и то же), но в двух разных точках. В нестохастическом случае достаточно иметь возможность одного обращения на каждом шаге.*

Число итераций (с точностью до множителя 2 в стохастическом случае) – это число обращений к такому оракулу. Цель статьи: обращаясь к оракулу на одном шаге (итерации) не более двух раз, так организовать итерационную процедуру, чтобы сгенерированная на основе опроса оракула последовательность $\{x^k\}$ с вероятностью не менее $1-\sigma$ удовлетворяла неравенству

$$f\left(\frac{1}{N}\sum_{k=1}^{N} x^k\right) - f_* \le \frac{1}{N}\sum_{k=1}^{N} f(x^k) - f_* \le \varepsilon$$

с как можно меньшим значением $N$.

**Замечание 3.** На самом деле, не очень важно, сколько раз разрешено обращаться к оракулу, важно только, что не менее двух раз [2]. Приведенные в разделах 4, 5 результаты легко переписываются, если вместо двух точек (на одной реализации) разрешается использовать $k \le n+1$ точек (на одной реализаций): грубо говоря, оценки числа итераций от желаемой точности улучшатся в $k$ раз $N(\varepsilon) \to N(\varepsilon)/k$ [2, 16]. Если же разрешается обращаться только один раз, то картина принципиально меняется [2]. В этом случае на данный момент имеется достаточно большой зазор (для детерминированных постановок задач) между нижними оценками и тем, что сейчас дают лучшие методы [2, 7, 17].

Изложим далее общую схему, позволяющую свести описанную выше постановку к постановке раздела 2. Тогда можно будет воспользоваться теоремой 1.

Пусть $e \in RS_p^n(1)$ ($\tilde{e} \in RB_p^n(1)$) – случайный вектор, равномерно распределенный на сфере (шаре) единичного радиуса в $l_p$ норме в $\mathbb{R}^n$ (далее ограничимся рассмотрением случаев: $p=1$, $p=2$, $p=\infty$). Сгладим (следуя [7]) исходную функцию с помощью локального усреднения по шару радиуса $\mu > 0$ ($\mu \le \mu_0$), который будет выбран позже,

$$f^\mu(x;\eta) = E_{\tilde{e}}\left[f(x + \mu\tilde{e};\eta)\right],$$
$$f^\mu(x) = E_{\tilde{e},\eta}\left[f(x + \mu\tilde{e};\eta)\right].$$

Заменим исходную задачу (1) задачей

(3) $$f^\mu(x) \to \min_{x \in S_n(1)}.$$

Легко проверить (см., например, [7, 18, 19, 20] для $p'=2$, в общем случае рассуждения в точности такие же), что если выполняется **условие**

4. $\left|f(x;\eta)-f(y;\eta)\right| \le M_{p'}(\eta)\|x-y\|_{p'}$, $M_{p'} = \sqrt{E_\eta\left[M_{p'}(\eta)^2\right]} < \infty$,

то

$$0 \le f^\mu(x;\eta) - f(x;\eta) \le M_{p'}(\eta)\mu,$$

$$0 \le f^\mu(x) - f(x) \le M_{p'}\mu.$$

Это условие обобщает условие 3 раздела 2, в частности, в условии 4 соответствует в условии 3.а раздела 2.

Если выполняется **условие**

5. $\left\|\nabla_x f(x;\eta) - \nabla_x f(y;\eta)\right\|_{q'} \le L_{p'}(\eta)\|x-y\|_{p'}$, $L_{p'} = \sqrt{E_\eta\left[L_{p'}(\eta)^2\right]} < \infty$,

то

$$0 \le f^\mu(x;\eta) - f(x;\eta) \le L_{p'}(\eta)\mu^2/2,$$

$$0 \le f^\mu(x) - f(x) \le L_{p'}\mu^2/2,$$

где $1/p' + 1/q' = 1$. Предположим, что (если не предполагается гладкости, то можно просто положить $L_{p'} = \infty$)

(4) $$\min\left\{M_{p'}\mu, L_{p'}\mu^2/2\right\} \le \varepsilon/2,$$

и с вероятностью не менее $1-\sigma$ удалось получить следующее неравенство (например, воспользовавшись каким-то образом для задачи (3) теоремой 1):

$$\frac{1}{N}\sum_{k=1}^N f^\mu(x^k) - \min_{x \in S_n(1)} f^\mu(x) \le \frac{\varepsilon}{2}.$$

Тогда с вероятностью $\ge 1-\sigma$:

$$f\left(\frac{1}{N}\sum_{k=1}^N x^k\right) - f_* \le \frac{1}{N}\sum_{k=1}^N f(x^k) - \min_{x \in S_n(1)} f(x) \le \frac{1}{N}\sum_{k=1}^N f^\mu(x^k) - \min_{x \in S_n(1)} f^\mu(x) + \frac{\varepsilon}{2} \le \varepsilon.$$

Таким образом, при условии (4) решение задачи (3) с точностью $\varepsilon/2$ является решением задачи (1) с точностью $\varepsilon$.

Сглаживание было введено для того, чтобы для сглаженной задачи с помощью описанного оракула можно было получить несмещенную оценку субградиента. К сожале-нию, без сглаживания не понятно, как это можно было бы сделать. Итак, введем (при $p = 2$, см., например, [19]) аналог стохастического субградиента

$$g^{\mu}(x;e,\eta) = \frac{\mathrm{Vol}(S_p^n(\mu))}{\mathrm{Vol}(B_p^n(\mu))}(f(x+\mu e;\eta) - f(x;\eta))\bar{e},$$

где $e$ – случайный вектор, равномерно распределенный на сфере радиуса 1 в $l_p$ норме (обозначим такую сферу через $S_p^n(1)$); $\mathrm{Vol}(B_p^n(\mu))$ – объем шара радиуса $\mu$ в $l_p$ норме, аналогично определяется $\mathrm{Vol}(S_p(\mu))$; $\bar{e} = \bar{e}(e)$ – вектор с $l_2$ нормой, равной 1, ортогональный поверхности $S_p^n(1)$ в точке $e$. Например, см. табл.1.

**Таблица 1**

| $p$ | Аналог стохастического субградиента | Выбор направления |
|---|---|---|
| 1 | $\dfrac{n}{\mu}(f(x+\mu e;\eta) - f(x;\eta))\begin{pmatrix} \mathrm{sign}\,e_1 \\ \ldots\ldots \\ \mathrm{sign}\,e_n \end{pmatrix}$ | $e \in RS_1^n(1)$ |
| 2 | $\dfrac{n}{\mu}(f(x+\mu e;\eta) - f(x;\eta))e$ | $e \in RS_2^n(1)$ |
| $\infty$ | $\dfrac{n}{\mu}(f(x+\mu e;\eta) - f(x;\eta))\breve{e}_{i(e)}$ (п.н.) $\breve{e}_{i(e)} = (\underbrace{0,\ldots,0,1,0,\ldots,0}_{i(e)})$, $i(e) = \arg\max\limits_{i=1,\ldots,n}|e_i|$ | $e \in RS_\infty^n(1)$ |

Основное свойство $g^{\mu}(x;e,\eta)$ заключается в том, что

$$E_{e,\eta}[g^{\mu}(x;e,\eta)] \equiv \nabla f^{\mu}(x).$$

Воспользовались векторным вариантом теоремы Стокса, подобно [20], см. Приложение 1 для $p = 1$.

Причем это свойство сохраняется и в случае, когда вместо "идеального" значения реализаций $f(x+\mu e;\eta)$ и $f(x;\eta)$ оракул выдает зашумленные

$$E_{e,\eta}[g_\delta^{\mu}(x;e,\eta)] \equiv \nabla f^{\mu}(x).$$

Чтобы можно было воспользоваться теоремой 1 для сглаженной задачи (3), необходимо оценить $\|g_\delta^{\mu}(x;e,\eta)\|_\infty$, где

$$g_\delta^\mu(x;e,\eta) = \frac{\text{Vol}(S_p^n(\mu))}{\text{Vol}(B_p^n(\mu))}\Big(f(x+\mu e;\eta) + \tilde{\delta}_{x+\mu e}(\eta) - \big(f(x;\eta) + \tilde{\delta}_x(\eta)\big)\Big)\overline{e}.$$

Из определения оракула следует, что при $p=1$ и условии 3 раздела 2

(5) $$\left\|g_\delta^\mu(x;e,\eta)\right\|_\infty \leq \left(M + \frac{2\delta}{\mu}\right)n.$$

При $p=2$ и $p=\infty$ оценка (5) получается хуже (см. табл.2 ниже).

Выберем согласно условию (4) $\mu = \varepsilon/(2M)$ и будем считать, что (условие на допустимый уровень шума)

$$\delta \leq \varepsilon/4.$$

Тогда условие (5) перепишется следующим образом:

$$\left\|g_\delta^\mu(x;e,\eta)\right\|_\infty \leq 2Mn.$$

Подобно алгоритму (2) опишем алгоритм решения задачи (3) для $p=1$. Положим $x_i^1 = 1/n$, $i=1,...,n$. Пусть $t=1,...,N-1$.

$$x_i^{t+1} = \frac{\exp\left(-\frac{1}{\beta_{t+1}}\sum_{k=1}^{t}\left[g_\delta^\mu(x^k;e^k,\eta^k)\right]_i\right)}{\sum_{l=1}^{n}\exp\left(-\frac{1}{\beta_{t+1}}\sum_{k=1}^{t}\left[g_\delta^\mu(x^k;e^k,\eta^k)\right]_l\right)}, \quad i=1,...,n, \quad \beta_t = \frac{2Mn\sqrt{t}}{\sqrt{\ln n}},$$

где $[z]_i$ – $i$-я координата вектора $z$.

**Теорема 2.** *Пусть имеется оракул из предположений 1, 2 с $\delta \leq \varepsilon/4$ и справедливы условия 1, 2, 3.а раздела 2, $p=1$, тогда для*

$$N = \left\lceil \frac{64M^2 n^2 \ln n}{\varepsilon^2} \right\rceil$$

*имеет место оценка*

$$E\left[f\left(\frac{1}{N}\sum_{k=1}^{N}x^k\right)\right] - f_* \leq \varepsilon.$$

*Если (дополнительно) справедливо условие 3 раздела 2, тогда для*

$$N = \left\lceil \frac{128M^2 n^2}{\varepsilon^2}\left(\ln n + 8\ln(\sigma^{-1})\right) \right\rceil$$

*с вероятностью не менее $1-\sigma$ имеет место оценка*

$$f\left(\frac{1}{N}\sum_{k=1}^{N}x^k\right) - f_* \leq \varepsilon.$$

**Доказательство.** Применим теорему 1 к функции $f^\mu(x)$ с

$$N = \left\lceil \frac{4(2Mn)^2}{(\varepsilon/2)^2} \ln n \right\rceil$$

для оценки скорости сходимости по математическому ожиданию и с

$$N = \left\lceil \frac{8(2Mn)^2}{(\varepsilon/2)^2}\left(\ln n + 8\ln(\sigma^{-1})\right) \right\rceil = \left\lceil \frac{128M^2n^2}{\varepsilon^2}\left(\ln n + 8\ln(\sigma^{-1})\right) \right\rceil$$

для оценки скорости сходимости с учетом вероятностей больших уклонений. В последнем случае было использовано неравенство $\left(\sqrt{a}+\sqrt{b}\right)^2 \leq 2a+2b$. □

Резюмируем полученные результаты в виде таблицы. При этом считаем выполненными условия 1, 2 раздела 2 и условие 4 (константы $M_1$, $M_2$, $M_\infty$ определяются в условии 4). Во второй строчке таблицы приведены математические ожидания числа итераций. Заметим при этом, что

$$M_1^2 \leq M_2^2 \leq nM_1^2, \; M_2^2 \leq M_\infty^2 \leq nM_2^2.$$

**Таблица 2**

| $p=1$ | $p=2$ | $p=\infty$ |
|---|---|---|
| $\mathrm{O}\left(\dfrac{M_1^2 n^2 \ln n}{\varepsilon^2}\right)$ | $\mathrm{O}\left(\dfrac{M_2^2 n^2 \ln n}{\varepsilon^2}\right)$ | $\mathrm{O}\left(\dfrac{M_\infty^2 n^2 \ln n}{\varepsilon^2}\right)$ |

Неизвестно, оптимальна ли выписанная в табл.2 оценка для $p=1$ при наложенных условиях на уровень шума $\delta \leq \varepsilon/4$. Однако имеется гипотеза, что полученная оценка оптимальна с точностью до мультипликативной константы при заданном уровне шума. В условиях отсутствия шума ($\delta = 0$) приведенная оценка (и тем более остальные для $p=2$, $p=\infty$), вообще говоря, не является оптимальной для негладких задач стохастической оптимизации. В [2] получена оценка (с помощью техники двойного сглаживания), которая позволяет сократить число итераций в оценке из второго столбца табл.2 ($p=2$) в $\sim n/\ln n$ раз (причем, по-видимому, логарифмический множитель тут можно убрать). Однако предложенный в [2] метод непрактичный (в отличие от предложенного в данном разделе метода), поскольку чрезвычайно чувствителен даже к очень небольшим шумам.

## 4. Модификация метода зеркального спуска для гладких задач стохастической оптимизации при спусках по случайному направлению

К сожалению, описанный в разделе 3 подход дает оценку в $n$ раз большую нижней оценки в гладком случае [2]. Поскольку интересны ситуации, в которых $n \gg 1$, то необходимо этот зазор как-то устранить. Естественно попытаться найти в рассуждениях раздела 3 наиболее грубое место и попробовать провести более точные рассуждения. К счастью, такое место всего одно – неравенство (5).

Считаем далее выполненными условия 1, 2 раздела 2 и условие 4 раздела 3 (с $p' = 2$).

Чтобы пояснить, в чем заключается грубость, рассмотрим для большей наглядности случай с $\delta = 0$. Тогда можно устремить $\mu \to 0+$ и получить

$$g^\mu(x;e,\eta) \to g(x;e,\eta) = \frac{\text{Vol}(S_p^n(1))}{\text{Vol}(B_p^n(1))} \langle \nabla_x f(x;\eta), e \rangle \overline{e}.$$

Аналогично разделу 3 имеем

$$E_{e,\eta}\left[g(x;e,\eta)\right] \equiv \nabla f(x).$$

Оценим сначала $E_{e,\eta}\left[\|g(x;e,\eta)\|_{\overline{q}}^2\right]$ ($2 \leq \overline{q} \leq \infty$ выбирается исходя из структуры множества, на котором происходит оптимизация,[4] см. [3]) при $p = 1$ (этот параметр отве-

---

[4] В разбираемом в статье случае, когда ограничение в виде симплекса, выбирают $\overline{q} = \infty$, см. раздел 2 (обоснование такому выбору имеется, например, в [21]). Как уже отмечалось, выбор $\overline{q}$ осуществляется исходя из структуры множества, на котором происходит оптимизация. Вместо используемого нами варианта метода зеркального спуска (МЗС) из раздела 2, "настроенного" на то, что оптимизация происходит на симплексе, можно использовать вариант, подходящий для любого другого выпуклого множества $Q$, в котором в прямом пространстве выбрана норма $l_{\overline{p}}$ ($1/\overline{p} + 1/\overline{q} = 1$) и определена неотрицательная сильно выпуклая (с константой $\geq 1$) относительно этой нормы функция $d(x)$, задающая "расстояние" Брэгмана

$$V(x,y) = d(x) - d(y) - \langle \nabla d(y), x - y \rangle.$$

Итоговая оценка ожидаемого числа итераций для соответствующего МЗС (в случае, когда на каждой итерации доступен несмещенный стохастический субградиент, математическое ожидание квадрата $l_{\overline{q}}$ нормы которого равномерно по $x$ ограничено числом $M_{\overline{p}}^2$) будет [3, чает за выбор способа рандомизации в табл.1 раздела 3) в категориях $\mathrm{O}(\ )$. Для этого заметим,[5] что случайный вектор $e \in RS_1^n(1)$ можно представить как $e = a/\|a\|_1$, где компоненты вектора – независимые лапласовские случайные величины, т.е. с плотностью $e^{-|y|}/2$. Согласно табл.1 имеем

$$E_{e,\eta}\left[\|g(x;e,\eta)\|_{\bar{q}}^2\right] = n^{2+2/\bar{q}} E_{e,\eta}\left[\langle\nabla_x f(x;\eta), e\rangle^2\right] = n^{2+2/\bar{q}} E_{e,\eta}\left[\frac{\langle\nabla_x f(x;\eta), a\rangle^2}{\|a\|_1^2}\right].$$

Далее воспользуемся тем, что $n \gg 1$. Тогда исходя из явления концентрации меры[6] [23, 24], имеем: $\|a\|_1^2$ – сконцентрирован (с хвостами вида $e^{-\sqrt{y}}$) около своего математического ожидания $c_1 n$; $\langle\nabla_x f(x;\eta), a\rangle^2$ – $e^{-\sqrt{y}}$ – сконцентрирован (при зафиксированном $\eta$) около своего математического ожидания $c_2 \|\nabla_x f(x;\eta)\|_2^2$.[7] В результате получается следующая оценка (отметим, что здесь и далее $M_2$ определяется условием 4 раздела 3)

(6) $$E_{e,\eta}\left[\|g(x;e,\eta)\|_{\bar{q}}^2\right] = \mathrm{O}\left(n^{1+2/\bar{q}} M_2^2\right) \text{ (при } p = 1\text{).}$$

В действительности можно показать, что и $\|g(x;e,\eta)\|_{\bar{q}}^2$ имеет $e^{-\sqrt{y}}$ – концентрацию около своего математического ожидания, если в условии 4 раздела 3 $M_2(\eta) \equiv M_2$.

---

21]: $\mathrm{O}\left(M_{\bar{p}}^2 V(x^*, x^1)/\varepsilon^2\right)$. Причем если $Q = B_{\bar{p}}(R)$, то оптимально (с точностью до мультипликативного множителя) выбирать:

$$d(x) = \frac{1}{2(a-1)}\|x\|_a^2, \ 1 \leq \bar{p} \leq a; \ d(x) = \frac{1}{2(\bar{p}-1)}\|x\|_{\bar{p}}^2, \ a \leq \bar{p} \leq 2; \ d(x) = \frac{1}{2}\|x\|_2^2, \ 2 \leq \bar{p} \leq \infty.$$

В разделе 2 использовался вариант МЗС с прокс-функцией $d(x) = \ln n + \sum_{k=1}^{n} x_k \ln x_k$, что также приводит к неулучшаемой (с точностью до числового множителя) оценке числа итераций [21].

[5] Приводимая далее схема рассуждений была сообщена Александром Содиным.

[6] Впрочем, можно приведенные ниже результаты получить и без тонких оценок плотности концентрации исходя из классических вариантов закона больших чисел, центральной предельной теоремы и их идемпотентных аналогов [22].

[7] Точные значения положительных констант $c_1$ и $c_2$ (аналогично $c_3$, $c_4$, $c_5$) здесь не интересны, сейчас важно только то, что они не зависят от $n$. Здесь и далее для большей наглядности предполагаем выполненным максимально сильное условие 3 раздела 2.

Еще более геометрически наглядные рассуждения, восходящие к Пуанкаре–Леви [23], связанные с концентрацией равномерной меры на евклидовой сфере, позволяют получить следующую оценку:

(7) $$E_{e,\eta}\left[\|g(x;e,\eta)\|_{\bar{q}}^2\right] = \mathrm{O}\left(n^{2/\bar{q}}\ln n M_2^2\right) \text{ (при } p=2\text{).}$$

Отличие в рассуждениях в том, что $e \in RS_2^n(1)$ стоит представлять как $e = a/\|a\|_2$, где $a \in \mathrm{N}(0, I_n)$, а $I_n$ – единичная матрица (на диагонали единицы, остальные элементы нули) размера $n \times n$. Тогда

$$E_{e,\eta}\left[\|g(x;e,\eta)\|_{\bar{q}}^2\right] = n^2 E_{e,\eta}\left[\langle \nabla_x f(x;\eta), e\rangle^2 \|e\|_{\bar{q}}^2\right] = n^2 E_{e,\eta}\left[\frac{\langle \nabla_x f(x;\eta), a\rangle^2 \|a\|_{\bar{q}}^2}{\|a\|_2^4}\right],$$

где $\|a\|_2^4 - e^{-\sqrt{y}}$ – сконцентрирован около своего математического ожидания $c_3 n^2$, $\|a\|_{\bar{q}}^2 - e^{-y}$ – сконцентрирован (экспоненциально сконцентрирован) около своего математического ожидания, которое оценивается сверху[8] величиной $c_4 n^{2/\bar{q}} \ln n$ [22], $\langle \nabla_x f(x;\eta), a\rangle^2$ – экспоненциально сконцентрирован (при зафиксированном $\eta$) около своего математического ожидания $c_5 \|\nabla_x f(x;\eta)\|_2^2$, которое оценивается сверху $c_5 M_2^2$, если в условии 4 раздела 3 $M_2(\eta) \equiv M_2$.

Наиболее же просто исследуется случай $p = \infty$. Основным здесь является следующее наблюдение: практически весь объем многомерного куба сосредоточен на его границе [23].[9] Таким образом, в предположении $n \gg 1$ с хорошей точностью мы можем заменить условие $e \in RS_\infty^n(1)$ условием $e \in RB_\infty^n(1)$. Последнее распределение тривиально исследуется [23]. Аналогично вышеизложенному

$$E_{e,\eta}\left[\|g(x;e,\eta)\|_{\bar{q}}^2\right] = n^2 E_{e,\eta}\left[\langle \nabla_x f(x;\eta), e\rangle^2 \|e\|_{\bar{q}}^2\right] = n^2 E_{e,\eta}\left[\langle \nabla_x f(x;\eta), e\rangle^2\right].$$

Таким образом,

---

[8] Эту оценку можно уточнить. В частности (бакалаврский диплом И.Н. Усмановой, МФТИ, 2015),

$$E_e\left[\|e\|_{\bar{q}}^2\right] \leq (\bar{q}-1) n^{2/\bar{q}-1}, \ E_e\left[\|e\|_\infty^2\right] \leq (4\ln n)/n, \ e \in RS_2^n(1).$$

Используем это далее, см. табл.4 и выкладки в разделе 5.

[9] Действительно, объем $n$-мерного куба со стороной 1 равен 1, а со стороной $1-\delta$ равен $(1-\delta)^n \ll 1$ – при достаточно больших $n$.

(8) $$E_{e,\eta}\left[\left\|g(x;e,\eta)\right\|_{\bar{q}}^2\right] = \mathrm{O}(n^2 M_2^2) \text{ (при } p = \infty\text{)}.$$

В действительности, можно показать, что и $\left\|g(x;e,\eta)\right\|_{\bar{q}}^2$ имеет экспоненциальную концентрацию (при зафиксированном $\eta$) около своего математического ожидания, если в условии 4 раздела 3 $M_2(\eta) \equiv M_2$.

Исходя из несмещенной оценки субградиента $g(x;e,\eta)$ можно построить алгоритм, аналогичный (2): в (2) заменяем $\nabla_x f(x;\eta)$ на $g(x;e,\eta)$. Подставляя в оценки (6) – (8) $\bar{q} = \infty$ (что соответствует рассматриваемой в данной статье оптимизации на симплексе [21]), получим итоговые оценки среднего числа итераций такого алгоритма (табл.3)

**Таблица 3**

| $p = 1$ | $p = 2$ | $p = \infty$ |
|---|---|---|
| $\mathrm{O}\left(\dfrac{nM_2^2 \ln n}{\varepsilon^2}\right)$ | $\mathrm{O}\left(\dfrac{M_2^2 \ln^2 n}{\varepsilon^2}\right)$ | $\mathrm{O}\left(\dfrac{n^2 M_2^2 \ln n}{\varepsilon^2}\right)$ |

Из табл.3 хорошо видно, какая рандомизация предпочтительнее: на евклидовой сфере ($p = 2$). Отсюда с учетом того, что $M_1^2 \le M_2^2 \le nM_1^2$, получаем оценку

$$\mathrm{O}\left(\frac{M_2^2 \ln^2 n}{\varepsilon^2}\right) \le \mathrm{O}\left(\frac{M_1^2 n \ln^2 n}{\varepsilon^2}\right),$$

которая с точностью до логарифмического множителя соответствует нижней оценке [2]. Однако если предположить, что $M_2^2 \ll nM_1^2$, то получается, что можно превзойти нижнюю оценку, т.е. быстрее достичь желаемой точности, чем предписано нижней оракульной оценкой [2, 7]. Но никакого противоречия здесь, конечно, нет, поскольку нижняя оценка была получена без всяких дополнительных предположений. Делая такое предположение ($M_2^2 \ll nM_1^2$), уже не в праве говорить об оценке [2] как о нижней оценке для этого нового класса.

Рассмотрены только три значения $p$ и только симплекс в качестве множества, на котором происходит оптимизация. Можно показать, что общий вывод сохранится при рассмотрении всевозможных $1 \le p \le \infty$ и всевозможных выпуклых множеств $Q$, на которых происходит оптимизация: *наиболее предпочтительная рандомизация* $e \in RS_2^n(1)$.

В действительности получен намного более общий результат (см. также [14]). Пусть рассматривается задача

$$f(x) = E_\eta\left[f(x;\eta)\right] \to \min_{x \in Q},$$

где $Q$ – выпуклое множество (необязательно ограниченное). Пусть в прямом пространстве выбрана $l_{\bar{p}}$ норма,[10] $1/\bar{p} + 1/\bar{q} = 1$. Введена соответствующая этой норме прокс-функция [3, 14, 21]. Пусть $R_{\bar{p}}^2$ – "расстояние" Брэгмана от точки старта до решения, посчитанное согласно этой прокс-функции [3, 14, 21]. Приводимая ниже табл.4 была ранее известна при $\bar{p} = \bar{q} = 2$ (при $\bar{q} = 2$ в сильно выпуклом случае можно убрать множитель $\ln n$).

**Таблица 4**

| Выполнены условия 1, 2 раздела 2 и условие[11] 4 раздела 3 | $f(x)$ – выпуклая функция | $f(x)$ – $\gamma_{\bar{p}}$-сильно выпуклая функция в $l_{\bar{p}}$ норме |
|---|---|---|
| $2 \leq \bar{q} \leq \ln n$ | $\mathrm{O}\left(\dfrac{\bar{q} M_2^2 R_{\bar{p}}^2 n^{2/\bar{q}}}{\varepsilon^2}\right)$ | $\mathrm{O}\left(\dfrac{\bar{q} M_2^2 n^{2/\bar{q}} \ln n}{\gamma_{\bar{p}} \varepsilon}\right)$ |
| $\ln n \leq \bar{q} \leq \infty$ | $\mathrm{O}\left(\dfrac{M_2^2 R_{\bar{p}}^2 \ln n}{\varepsilon^2}\right)$ | $\mathrm{O}\left(\dfrac{M_2^2 \ln^2 n}{\gamma_{\bar{p}} \varepsilon}\right)$ |

Все сказанное выше относилось не к безградиентным методам, а к методам спуска по случайному направлению и притом в гладком случае. Однако нижние оценки тут с точностью до логарифмических множителей одинаковы. Выше было показано, как можно для спусков по случайному направлению в гладком случае приблизиться, а в определенных ситуациях и превзойти нижнюю оценку. Естественно возникает желание перенести предложенный здесь оптимальный метод и на безградиентные методы так, чтобы сохранить полученную оценку. При этом необходимо определить уровень допустимого шума, при котором это возможно. Собственно, этому и посвящен следующий раздел.

Сейчас же остановимся на одном интересном обстоятельстве, выявленном в разделах 3, 4. Получается довольно неожиданная ситуация: оказывается, имеет место сильная зависимость скорости сходимости метода от того, какой способ рандомизации (а по сути сглаживания) выбирать. Причем, как это видно из табл.3, разница очень существенная. К сожалению, в своем желании сохранить несмещенность оценки субградиента мы "перегнули палку" в случае $p=1$ и особенно $p=\infty$. Несмещенность в этих случаях досталась дорогой ценой – большой оценкой дисперсии соответствующих оценок. Собственно, при предельном переходе $\mu \to 0+$ была унаследована большая дисперсия, что и наблюдали в табл.3. Естественно, в этой связи задаться вопросом: а может быть рандомизация $e \in RS_2^n(1)$ оптимальна только в классе несмещенных оценок? А если допускать смещение (bias), то, возможно, можно будет добиться лучшего, как, скажем, в случае оптимальных оценок в математической статистике [26] (см. также [20])? Оказывается, что если допускать смещение, рандомизация $e \in RS_2^n(1)$ по-прежнему будет оптимальной (с точностью до логарифмического множителя). Чтобы это пояснить, продолжим рассмотрение гладкого случая с возможностью получения на каждом шаге (итерации) от оракула незашумленной производной по указанному направлению. Рассмотрим более общую схему (см., например, [2]). Пусть $Z$ – случайный вектор с корреляционной матрицей $E_Z\left[ZZ^T\right]=I_n$. Тогда

$$g(x;Z,\eta) = \langle \nabla_x f(x;\eta), Z \rangle Z = ZZ^T \nabla_x f(x;\eta).$$

---

[10] $\bar{p} \in [1, 2]$ – другие значения, как правило, не интересны [3, 14, 21].

[11] С выполнением условия 4 есть нюанс, когда $Q$ не ограничено [14]. Однако можно искусственно компактифицировать $Q$ исходя из того, что по ходу итерационного процесса "расстояние" от текущей точки до решения может быть оценено сверху "расстоянием" от точки старта до решения, умноженным на некоторую константу (см. доказательство теоремы 4 в [25]).

Очевидно, что

$$E_{Z,\eta}\left[g(x;Z,\eta)\right] = \nabla f(x).$$

Оказывается, можно улучшить оценку, соответствующую $p=1$, выбирая в этом подходе случайный вектор $Z$ так, чтобы каждая компонента принимала независимо и равновероятно одно из двух значений $1$, $-1$ (равномерное распределение на Хэмминговском кубе), см. табл.1. Тогда[12] [2] (см. табл.3)

$$O\left(\frac{nM_2^2 \ln n}{\varepsilon^2}\right) \to O\left(\frac{M_2^2 \ln n}{\varepsilon^2}\right),$$

что улучшает приведенную ранее оценку с рандомизацией на евклидовой сфере на логарифмический множитель.

---

[12] Основная выкладка, поясняющая формулу, достаточно тривиальна

$$E_{Z,\eta}\left[\|g(x;Z,\eta)\|_\infty^2\right] = E_{Z,\eta}\left[\|\langle \nabla_x f(x;\eta), Z\rangle Z\|_\infty^2\right] = E_{Z,\eta}\left[|\langle \nabla_x f(x;\eta), Z\rangle|^2\right] =$$

$$= E_{Z,\eta}\left[\langle \nabla_x f(x;\eta), Z\rangle^2\right] = E_\eta\left[\nabla_x f(x;\eta)^T \underbrace{E_Z\left[ZZ^T\right]}_{I_n} \nabla_x f(x;\eta)\right] \leq M_2^2.$$

Если выбрать $Z \in \sqrt{n}\,\mathrm{N}(0, I_n)$, то получим в точности те же самые оценки, что получали ранее с рандомизацией на евклидовой сфере.

Если выбрать[13] (см. табл.1)

$$Z = \sqrt{n}\,(\,\underbrace{0,...,0,1}_{i},0,...,0\,),$$

где случайная величина $i$ независимо и равновероятно принимает значения $1,...,n$, т.е.

$$P(i = k) = \frac{1}{n},\ k = 1,...,n,$$

и считать, что в прямом пространстве выбрана (в связи со свойствами множества $Q$) норма $l_2$, то такая покомпонентная рандомизация приводит к аналогичным оценкам, даваемым рандомизацией на евклидовой сфере при $\bar{p} = \bar{q} = 2$, что отражено в табл.5 [2].

**Таблица 5**

| $\bar{p} = \bar{q} = 2$ | $f(x)$ – выпуклая функция | $f(x)$ – $\gamma_2$-сильно выпуклая функция в $l_2$ норме |
|---|---|---|
| Рандомизация на евклидовой сфере | $\mathrm{O}\!\left(\dfrac{M_2^2 R_2^2 n}{\varepsilon^2}\right)$ | $\mathrm{O}\!\left(\dfrac{M_2^2 n}{\gamma_2 \varepsilon}\right)$ |
| Покомпонентная рандомизация | $\mathrm{O}\!\left(\dfrac{M_2^2 R_2^2 n}{\varepsilon^2}\right)$ | $\mathrm{O}\!\left(\dfrac{M_2^2 n}{\gamma_2 \varepsilon}\right)$ |

Таким образом, видно, что в случае рассмотрения методов спуска по случайному направлению (покомпонентных методов) вполне можно рассчитывать на альтернативный способ получения оптимальных методов (оценок). Причем в последнем случае (покомпонентной рандомизации), можно существенно выиграть в стоимости одной итерации. Ранее такая задача в этой статье не ставилась. Было стремление минимизировать число обращений к оракулу (за значением функции, за производной по направлению), гарантирующих достижение заданной точности по функции. Если же минимизировать общую вычислительную сложность (число арифметических операций), то покомпонентные методы для большого класса важных в приложениях задач позволяют эффективно организовывать пересчет компонент градиента, т.е. не рассчитывать их каждый раз заново, что позволяет серьезно сэкономить в общих трудозатратах по сравнению с рандомизацией на евклидо-вой сфере (см., например, [14, 27, 28, 29]). Отметим также, что в последней строчке табл.5 константу $M_2^2$ можно оценивать как среднее значение по направлениям координатных осей, в то время как в предпоследней строчке таблицы $M_2^2$ оценивается по худшему направлению [29].

---

[13] Именно с дискретных аналогов такого подхода и начиналось изучение безградиентных методов [4, 5, 19].

Сказанное выше относилось к методам спуска по случайному направлению. Оказывается [2], что эти результаты можно перенести и на безградиентные методы. Для этого вводится аналог $g(x;Z,\eta)$:

$$g^\tau(x;Z,\eta) = \frac{f(x+\tau Z;\eta) - f(x;\eta)}{\tau} Z,$$

аналогично разделу 3 можно ввести и шумы $g_\delta^\tau(x;Z,\eta)$. К сожалению, даже при $\delta = 0$ не получается несмещенность, т.е. не выполняется условие 2 раздела 2, необходимое для справедливости теоремы 1 раздела 2, которой пользуемся. Но у теоремы 1 есть обобщение (см., например, [30]) не только на произвольные выпуклые множества $Q$ (что ранее уже неявно использовалось при заполнении табл.4, 5), но и на случай, когда условие 2 выполняется неточно (это как раз нужный случай). Исходя из такого обобщения [30] можно перенести (без изменения) выписанные оценки (при условии достаточной малости $\tau$ и $\delta = 0$) на безградиентные методы [2], причем рассуждения [2] можно обобщить и на случай $\delta > 0$, контролируя уровень шума (не будем в этой статье приводить соответствующие выкладки). Более того, отмеченное обобщение (из [30]) теоремы 1 позволяет не делать никаких ограничений (типа предположения 1 раздела 3) на шум, кроме должной малости уровня шума. Сами оценки (числа итераций) при этом удается сохранить, но за счет ужесточения требований к уровню шума. Схематично детали такого обобщения описаны в [14, 15].

## 5. Перенесение результатов раздела 4 на безградиентные методы

Рассмотрим

$$g_\delta^\mu(x;e,\eta) = \frac{n}{\mu}\Big(f(x+\mu e;\eta) + \tilde{\delta}_{x+\mu e}(\eta) - \big(f(x;\eta) + \tilde{\delta}_x(\eta)\big)\Big)e,$$

где $e \in RS_2^n(1)$. Поскольку (см. раздел 3)

$$E_{e,\eta}\big[g_\delta^\mu(x;e,\eta)\big] \equiv \nabla f^\mu(x),$$

то для возможности использования теоремы 1 и схемы раздела 3 нужно аккуратно ограничить сверху (в случае $Q = S_n(1)$ имеем $\bar{q} = \infty$) $E_{e,\eta}\big[\|g_\delta^\mu(x;e,\eta)\|_{\bar{q}}^2\big]$. Далее сконцентриру-емся именно на этой задаче. Здесь не будем бороться за то, чтобы получить оценки вероятностей больших уклонений.

Рассмотрим гладкий случай, в данном случае это подразумевает, что дополнительно к условиям 1, 2 раздела 2 и условию 4 раздела 3 имеет место условие 5 раздела 3 (в обоих условиях предполагается, что $p' = 2$).

Из определения $g_\delta^\mu(x; e, \eta)$ и предположения 1 раздела 2 имеем

$$E_{e,\eta}\left[\left\|g_\delta^\mu(x; e, \eta)\right\|_{\bar{q}}^2\right] = \frac{n^2}{\mu^2} E_{e,\eta}\left[\left(f(x + \mu e; \eta) - f(x; \eta) + \left(\tilde{\delta}_{x+\mu e}(\eta) - \tilde{\delta}_x(\eta)\right)\right)^2 \|e\|_{\bar{q}}^2\right] =$$

$$= \frac{n^2}{\mu^2} E_{e,\eta}\left[\left(\left(f(x + \mu e; \eta) - f(x; \eta) - \mu\langle\nabla_x f(x; \eta), e\rangle\right) + \mu\langle\nabla_x f(x; \eta), e\rangle + \right.\right.$$

(9)
$$\left.\left. + \left(\tilde{\delta}_{x+\mu e}(\eta) - \tilde{\delta}_x(\eta)\right)\right)^2 \|e\|_{\bar{q}}^2\right].$$

Поскольку

(10)
$$\left|f(x + \mu e; \eta) - f(x; \eta) - \mu\langle\nabla_x f(x; \eta), e\rangle\right| \leq L_2(\eta)\mu^2/2,$$

$$\left|\tilde{\delta}_{x+\mu e}(\eta) - \tilde{\delta}_x(\eta)\right| \leq 2\delta,$$

$$(a + b + c)^2 \leq 3a^2 + 3b^2 + 3c^2,$$

то

$$E_{e,\eta}\left[\left\|g_\delta^\mu(x; e, \eta)\right\|_{\bar{q}}^2\right] \leq \frac{3}{4}n^2 L_2^2\mu^2 E_e\left[\|e\|_{\bar{q}}^2\right] + 3n^2 E_{e,\eta}\left[\langle\nabla_x f(x; \eta), e\rangle^2 \|e\|_{\bar{q}}^2\right] + 12\frac{n^2\delta^2}{\mu^2}E_e\left[\|e\|_{\bar{q}}^2\right].$$

Наиболее интересные ситуации - это $\bar{q} = 2$ и $\bar{q} = \infty$:

$$E_{e,\eta}\left[\left\|g_\delta^\mu(x; e, \eta)\right\|_{\bar{q}}^2\right] \leq 3nM_2^2 + \frac{3}{4}n^2 L_2^2\mu^2 + 12\frac{n^2\delta^2}{\mu^2} \text{ (при } \bar{q} = 2\text{)},$$

$$E_{e,\eta}\left[\left\|g_\delta^\mu(x; e, \eta)\right\|_{\bar{q}}^2\right] \leq 4\ln n M_2^2 + 3n\ln n L_2^2\mu^2 + 48\frac{n\ln n\delta^2}{\mu^2} \text{ (при } \bar{q} = \infty\text{)}.$$

Выберем $\mu$ согласно условию (4) раздела 3 $L_2\mu^2/2 \leq \varepsilon/2$, т.е. $\mu \leq \sqrt{\varepsilon/L_2}$. Следующее условие на $\mu$ и на допустимый уровень шума $\delta$ получим исходя из желания обеспечить выполнение неравенства (константа 5 здесь выбрана для определенности)

$$E_{e,\eta}\left[\left\|g_\delta^\mu(x; e, \eta)\right\|_{\bar{q}}^2\right] \leq 5nM_2^2 \text{ (при } \bar{q} = 2\text{)},$$

$$E_{e,\eta}\left[\left\|g_\delta^\mu(x; e, \eta)\right\|_{\bar{q}}^2\right] \leq 5\ln n M_2^2 \text{ (при } \bar{q} = \infty\text{)}.$$

Отсюда можно получить

$$\mu = \min\left\{\max\left\{\frac{\varepsilon}{2M_2}, \sqrt{\frac{\varepsilon}{L_2}}\right\}, \frac{M_2}{L_2}\sqrt{\frac{4}{3n}}\right\}, \quad \delta \le \frac{M_2\mu}{\sqrt{12n}} \text{ (при } \bar{q}=2\text{)},$$

$$\mu = \min\left\{\max\left\{\frac{\varepsilon}{2M_2}, \sqrt{\frac{\varepsilon}{L_2}}\right\}, \frac{M_2}{L_2}\sqrt{\frac{1}{6n}}\right\}, \quad \delta \le \frac{M_2\mu}{\sqrt{96n}} \text{ (при } \bar{q}=\infty\text{)}.$$

Подобно алгоритму (2) раздела 2 опишем оптимальный алгоритм (см. также раздел 3, только в п. 3 используется другая рандомизация) для задачи (1) и оракула из предположения 1 раздела 2. Положим $x_i^1 = 1/n$, $i=1,...,n$. Пусть $t=1,...,N-1$.

$$x_i^{t+1} = \frac{\exp\left(-\frac{1}{\beta_{t+1}}\sum_{k=1}^{t}\left[g_\delta^\mu\left(x^k; e^k, \eta^k\right)\right]_i\right)}{\sum_{l=1}^{n}\exp\left(-\frac{1}{\beta_{t+1}}\sum_{k=1}^{t}\left[g_\delta^\mu\left(x^k; e^k, \eta^k\right)\right]_l\right)}, \quad i=1,...,n, \quad \beta_t = M_2\sqrt{5t},$$

где $[z]_i$ – $i$-я координата вектора $z$.

**Теорема 3.** *Пусть имеется оракулом из предположений 1, 2 с $\delta \le M_2\mu/\sqrt{96n}$ и справедливы условия 1, 2, раздела 2 и условия 4, 5 раздела 3 (в которых $p'=2$). Тогда для задачи (1) и описанного выше алгоритма при*

$$N = \left\lceil \frac{80M_2^2 \ln^2 n}{\varepsilon^2} \right\rceil$$

*имеет место оценка*

$$E\left[f\left(\frac{1}{N}\sum_{k=1}^{N}x^k\right)\right] - f_* \le \varepsilon.$$

**Доказательство.** Применим теорему 1 (с учетом выкладок раздела 3) с

$$N = \left\lceil \frac{4 \cdot 5 M_2^2 \ln n}{(\varepsilon/2)^2} \ln n \right\rceil$$

к функции $f^\mu(x)$. □

Согласно [2] эта оценка оптимальна с точностью до мультипликативного множителя. Более того, подобно разделу 4 можно заметить, что в определенных случаях полученная оценка будет лучше нижней (оптимальной) оценки [2].

При дополнительных условиях (уточняющих условия 4, 5) здесь так же, как и в разделе 3, можно получить оценки вероятностей больших уклонений, однако не будем здесь приводить соответствующие оценки.

## 6. Заключение

В работе предложены эффективные методы нулевого порядка (также говорят прямые методы или безградиентные методы) для задач стохастической выпуклой оптимизации на симплексе и более общих выпуклых множествах с хорошей проксимальной структурой. Методы строились на базе обычного зеркального спуска для задач стохастической оптимизации. Вместо стохастического градиента в алгоритм зеркального спуска подставлялась специальная конечная разность, аппроксимирующая стохастический градиент. При правильном пересчете размера шага, получался эффективный метод, работающий по известным нижним оценкам, и даже их немного улучшающий при определенных условиях.

Все полученные результаты, кроме третьего столбца табл.4 (здесь нам известен только результат для $\bar{p} = \bar{q} = 2$), переносятся на онлайн постановки [13]. Подробности будут изложены в следующей статье.

Оригинальность результатов работы обеспечивается за счет рассмотрения неточного оракула, выдающего на каждой итерации в двух разных точках зашумленные значение оптимизируемой функции на одной и той же реализации. Наличие указанных (дополнительных) шумов "моделирует" практическую неустойчивость конечного дифференцирования, положенного в основу практически всех безградиентных методов. В частности, наличие таких шумов "в первом приближении" может моделировать конечность длины мантиссы.

В отличие от большинства других работ в данной работе также прорабатывался вопрос оптимального сочетания способа рандомизации при конструировании "дискретного стохастического градиента", использующегося в методах вместо недоступного настоящего стохастического градиента, с выбором прокс-структуры, определяемой геометрией выпуклого множества, на котором происходит оптимизация. В частности, подробно рассматривался, пожалуй, наиболее интересный пример такого множества (после евклидова шара) – симплекс. Ответ оказался достаточно универсальным: для гладких задач оптимальная рандомизация (в независимости от структуры множества) – рандомизация на евклидовой сфере. Для негладких задач в случае достаточно больших шумов, этот ответ уже может быть не верен. Соответствующий пример разбирался в разделе 3.

## Литература

1. *Нестеров Ю.Е.* Алгоритмические модели человеческого поведения. Выступление на Математическом кружке. Москва, МФТИ & МЦНМО, 14 сентября 2012 г. http://www.mathnet.ru/php/seminars.phtml?option_lang=rus&presentid=6990

2. *Duchi J.C., Jordan M.I., Wainwright M.J., Wibisono A.* Optimal rates for zero-order convex optimization: the power of two function evaluations // IEEE Transaction of Information. 2015. V. 61. № 5. P. 2788–2806.
http://www.eecs.berkeley.edu/~wainwrig/Papers/DucZero15.pdf

3. *Agarwal A., Bartlett P.L., Ravikumar P., Wainwright M.J.* Information-theoretic lower bounds on the oracle complexity of stochastic convex optimization // IEEE Transaction of Information. 2012. V. 58. № 5. P. 3235–3249. arXiv:1009.0571

4. *Kiefer J., Wolfowitz J.* Statistical estimation on the maximum of a regression function // Ann. Math. Statist. 1952. V. 23. P. 462–466.

5. *Поляк Б.Т.* Введение в оптимизацию. М.: Наука, 1983.

6. *Граничин О.Н., Поляк Б.Т.* Рандомизированные алгоритмы оценивания и оптимизации при почти произвольных помехах. М.: Наука, 2003.

7. *Немировский А.С., Юдин Д.Б.* Сложность задач и эффективность методов оптимизации. М.: Наука, 1979.

8. *Konečný J., Richárik P.* Simple complexity analysis of simplified direct search // e-print, 2014. arXiv:1410.0390

9. *Shapiro A., Dentcheva D., Ruszczynski A.* Lecture on stochastic programming. Modeling and theory. MPS-SIAM series on Optimization, 2014.

10. *Nesterov Y.* Primal-dual subgradient methods for convex problems // Math. Program. Ser. B. 2009. V. 120(1). P. 261–283.

11. *Юдицкий А.Б., Назин А.В., Цыбаков А.Б., Ваятис Н.* Рекуррентное агрегирование оценок методом зеркального спуска с усреднением // Проблемы передачи информации. 2005. Т. 41:4. С. 78–96.

12. *Nemirovski A., Juditsky A., Lan G., Shapiro A.* Stochastic approximation approach to stochastic programming // SIAM Journal on Optimization. 2009. V. 19. № 4. P. 1574–1609.

13. *Гасников А.В., Нестеров Ю.Е., Спокойный В.Г.* Об эффективности одного метода рандомизации зеркального спуска в задачах онлайн оптимизации // ЖВМ и МФ. Т. 55. № 4. 2015. С. 55–71.

14. *Гасников А.В., Двуреченский П.Е., Нестеров Ю.Е.* Стохастические градиентные методы с неточным оракулом // Труды МФТИ. 2016. Т. 8. № 1. arxiv:1411.4218

15. *Bogolubsky L., Dvurechensky P., Gasnikov A., Gusev G., Nesterov Yu., Raigorodskii A., Tikhonov A., Zhukovskii M.* Learning supervised PageRank with gradient-based and gradient-free optimization methods // e-print, 2016. arXiv:1603.00717

16. *Гасников А.В., Двуреченский П.Е., Камзолов Д.И.* Градиентные и прямые методы с неточным оракулом для задач стохастической оптимизации // Динамика систем и процессы управления. Труды Международной конференции, посвящено 90-летию со дня рождения академика Н.Н. Красовского. Екатеринбург, 15 – 20 сентября 2014. Издательство: Институт математики и механики УрО РАН им. Н.Н. Красовского (Екатеринбург), 2015. С. 111–117. arXiv:1502.06259

17. *Belloni A., Liang T., Narayanan H., Rakhlin A.* Escaping the Local Minima via Simulated Annealing: Optimization of Approximately Convex Functions // JMLR. Workshop and conference proceedings. 2015. V. 40. P. 1 – 26. arXiv:1501.07242

18. *Nesterov Yu.* Random gradient-free minimization of convex functions // CORE Discussion Paper 2011/1. 2011.

19. *Spall J.C.* Introduction to stochastic search and optimization: estimation, simulation and control. Wiley, 2003.

20. *Bubeck S., Cesa-Bianchi N.* Regret analysis of stochastic and nonstochastic multi-armed bandit problems // Foundation and Trends in Machine Learning. 2012. V. 5. № 1. P. 1–122. arXiv:1204.5721

21. *Nemirovski A.* Lectures on modern convex optimization analysis, algorithms, and engineering applications. Philadelphia: SIAM, 2013.
http://www2.isye.gatech.edu/~nemirovs/Lect_ModConvOpt.pdf

22. *Лидбеттер М., Линдгрен Г., Ротсен Х.* Экстремумы случайных последовательностей и процессов. М.: Мир, 1989.

23. *Ledoux M.* Concentration of measure phenomenon. Providence, RI, Amer. Math. Soc., 2001 (Math. Surveys Monogr. V. 89).

24. *Boucheron S., Lugoshi G., Massart P.* Concentration inequalities: A nonasymptotic theory of independence. Oxford University Press, 2013.

25. *Гасников А.В., Двуреченский П.Е., Дорн Ю.В., Максимов Ю.В.* Численные методы поиска равновесного распределения потоков в модели Бэкмана и модели стабильной динамики // Математическое моделирование. 2016. Т. 28. (принята к печати) arXiv:1506.00293

26. *Ибрагимов И.А., Хасьминский Р.З.* Асимптотическая теория оценивания. М.: Наука, 1977.

27. *Wright S.J.* Coordinate descent algorithms // Optimization online, 2015.

http://www.optimization-online.org/DB_FILE/2014/12/4679.pdf

28. *Anikin A., Dvurechensky P., Gasnikov A., Golov A., Gornov A., Maximov Yu., Mendel M., Spokoiny V.* Modern efficient numerical approaches to regularized regression problems in application to traffic demands matrix calculation from link loads // Proceedings of International conference ITAS-2015. Russia, Sochi, September, 2015. arXiv:1508.00858

29. *Гасников А.В., Двуреченский П.Е., Усманова И.Н.* О нетривиальности быстрых (ускоренных) рандомизированных методов // Труды МФТИ. 2016. Т. 8. № 2. arXiv:1508.02182

30. *Juditsky A., Nemirovski A.* First order methods for nonsmooth convex large-scale optimization, I, II. In: Optimization for Machine Learning. Eds. S. Sra, S. Nowozin, S. Wright. MIT Press, 2012.

## Приложение 1

**Лемма.** $E_{e,\eta}\left[g^{\mu}(x;e,\eta)\right] \equiv \nabla f^{\mu}(x)$.

**Доказательство.** Распишем левую часть выражения в условиях леммы

$$E_{e,\eta}\left[g^{\mu}(x;e,\eta)\right] = \frac{n}{\mu} E_{e,\eta}\left[\left(f(x+\mu e;\eta) - f(x;\eta)\right)\hat{e}\right] =$$

$$= \frac{n}{\mu} E_{e,\eta}\left[f(x+\mu e;\eta)\hat{e}\right] - \frac{n}{\mu} E_{e,\eta}\left[f(x;\eta)\hat{e}\right] = \frac{n}{\mu} E_{e,\eta}\left[f(x+\mu e;\eta)\hat{e}\right] =$$

$$= \frac{n\sqrt{n}}{\mu \operatorname{Vol}(S_1^n(\mu))} \int_{S_1^n(\mu)} E_{\eta}\left[f(x+\mu e;\eta)\right]\frac{\hat{e}}{\sqrt{n}} d\sigma(e) = \frac{n\sqrt{n}}{\mu^2 \operatorname{Vol}(S_1^n(\mu))} \int_{B_1^n(\mu)} \nabla f(x+v)\,dv =$$

$$= \frac{n\sqrt{n}\,\nabla E_{\tilde{e}}\left[f(x+\mu\tilde{e})\right]\operatorname{Vol}(B_1^n(\mu))}{\mu \operatorname{Vol}(S_1^n(\mu))},$$

где

$$\hat{e} = \left(\operatorname{sign} e_1, \ldots, \operatorname{sign} e_n\right)^T, \ e \in S_1^n(\mu), \ \tilde{e} \in B_1^n(\mu).$$

Заметим, что

$$\left\{\operatorname{Vol}(S_1^n(\mu)) = 2^n \frac{\sqrt{n}}{(n-1)!}\mu^{n-1}, \operatorname{Vol}(B_1^n(\mu)) = 2^n \frac{\mu^n}{n!}\right\} \Rightarrow \frac{\operatorname{Vol}(B_1^n(\mu))}{\operatorname{Vol}(S_1^n(\mu))} = \frac{\mu}{n\sqrt{n}}.$$

Таким образом, получаем

$$\frac{n\sqrt{n}\,\nabla E_{\tilde{e}}\left[f(x+\mu\tilde{e})\right]\operatorname{Vol}(B_1^n(\mu))}{\mu \operatorname{Vol}(S_1^n(\mu))} = \nabla E_{\tilde{e}}\left[f(x+\mu\tilde{e})\right] = \nabla f^{\mu}(x),$$

где $\tilde{e} \in B_1^n(\mu)$. □